\documentclass[a4paper,10pt]{article}

\usepackage{amsmath,amsfonts,amssymb,amsthm,epsfig}

\def\H{\mathbf{H}}

\def\gaga{\mathbf{U}_{n}^-}

\def\uasn{\mathbf{U}_v^-(\widehat{\mathfrak{sl}}_n)}

\def\C{\mathbb{C}}
\def\Z{\mathbb{Z}}
\def\N{\mathbb{N}}
\def\A{\mathbb{A}}
\def\K{\mathbb{K}}
\def\S{{\mathbb{S}}}

\def\F{\mathbb{F}}
\def\vac{|0\rangle}
\def\qed{$\hfill \square$}

\def\m{\mathbf{m}}
\def\n{\mathbf{n}}
\def\f{\mathbf{f}}
\def\Si{\mathfrak{S}}

\def\AH{\widehat{\mathbf{H}}}

\def\U{\mathbf{U}}

\def\b{\mathbf{b}}

\def\c{\mathbf{c}}

\def\l{\mathbf{l}}

\newtheorem*{theo}{Theorem}
\newtheorem*{prop}{Proposition}
\newtheorem*{lem}{Lemma}
\newtheorem*{cor}{Corollary}
\newtheorem*{conj}{Conjecture}
\newtheorem{lemma}{Lemma}
\numberwithin{equation}{section}

\begin{document}
\title{On the center of affine Hecke algebras of type A}
\author{Olivier Schiffmann}
\date{}
\maketitle
\noindent

\paragraph{0.1 Introduction.} Let $G$ be a simple complex algebraic group. Let
$W$ be its Weyl group and $\widehat{W}$ the associated extended affine Weyl
group. Let $\widehat{\H}$ be the Iwahori-Hecke algebra of $\widehat{W}$. It
is well-known that $\widehat{\H}$ admits two presentations~: the
\textit{Coxeter} presentation which arises naturally when $\widehat{\H}$ is
realized as the convolution algebra $L(\widehat{G},I)$ of compactly supported
functions on a p-adic group
$\widehat{G}=G(\overline{\mathbb{Q}_p})$ which are bi-invariant under action of
the Iwahori subgroup $I$ (see \cite{IM}), and the \textit{Bernstein}
presentation, which arises when $\widehat{\H}$ is realized in the
$G^\vee\times \C^*$-equivariant K-theory
of the Steinberg variety associated to $G^\vee$ where $G^\vee$ is the
Langlands dual group (see \cite{Gi}). The interplay
between these two presentations is central in the Deligne-Langlands
correspondence for finite-dimensional irreducible representations of
$\widehat{\H}$.\\
\hbox to1em{\hfill}The center $Z(\widehat{\H})$ is easily described in the
K-theoretic picture : it
is spanned by the classes of the trivial (equivariant) bundles on $Z$. A
geometric construction of this center in the convolution algebra presentation
is given by Gaitsgory, \cite{Ga}. This is in turn inspired by work of
Beilinson, and Haines, Kottwitz and Rapoport in the framework of Shimura
varieties, see \cite{Hai1},\cite{Hai2}.
\paragraph{}In this paper we give an explicit expression for the central
elements of $\widehat{\H}$ in the Coxeter presentation when $G=GL(r)$
(Theorem 2.5). This expression generalizes those obtained by Haines in the
minuscule case, \cite{Hai2} and is in some sense more explicit than \cite{Ga}.
More generally, we obtain expressions for central
elements in the ``parabolic spherical'' Hecke algebras
$L(\widehat{G},P)$ where $P\supseteq I$ is a parahoric subgroup. In
particular, taking $P=K$ to be a maximal compact open subgroup recovers
Lusztig's
description \cite{L1} of the Satake isomorphism between $Z(\widehat{\H})$ and
the spherical algebra $\widehat{\H}_{sph}$ (in the case $G=GL(r)$).
\paragraph{}Our method is based on the Hall algebra of a cyclic quiver, on
Uglov's higher-level Fock spaces and on the theory of canonical bases of
Kashiwara and Lusztig. Namely,
we use Ginzburg and Vasserot's geometric description of quantum affine
Schur-Weyl duality to construct an embedding of (half of) the center
$Z(\widehat{\H})$ in the center of the Hall algebra $\U^-_n$ of the quiver
$\widetilde{A}_{n-1}$ for $n \geq r$ (see \cite{S}). This
embedding is compatible with the canonical bases of $\widehat{\H}$ and
$\U^-_n$.
To describe the center of $\U^-_n$ we then consider the action on the Fock
spaces $\Lambda^\infty_{\mathbf{s}_l}$ recently introduced by Uglov \cite{U},
and use the
fact that this action is again compatible with the canonical bases.

\paragraph{}Finally, we give a simple alternate description of the center of
$\U^-_n$ in terms of a certain desingularization of orbit closures of
representations of the quiver $\widetilde{A}_{n-1}$, introduced by Varagnolo
and Vasserot \cite{VV}. This can be seen as a cyclic analogue of the
desingularization of orbit closures recently obtained by Reineke \cite{Re} for
finite-type simply laced Dynkin quivers.

\paragraph{}We note that the Fock spaces and their canonical bases
appear to be a very fundamental object in type A representation theory : they
describe Grothendieck groups and decomposition numbers of Hecke
algebras of type A or B (or more generally cyclotomic Hecke algebras) at
roots of unity (see
\cite{LLT},\cite{Ariki}, \cite{AM}, \cite{Gro}), and modular representations of
symmetric groups (see \cite{Dip}, \cite{Jam}, \cite{Gro}).

\paragraph{0.2 Notations.} Set $\S=\C[v]$, $\A=\C[v,v^{-1}]$ and $\K=\C(v)$. 
We define a $\C$-linear ring involution $u \mapsto \overline{u}$ on $\A$ by 
setting $\overline{v}=v^{-1}$.
Let $\mathbb{F}$ be a finite field with $q^2$ elements.
Let $\Si_r$ denote the
symmetric group on $r$ elements and let $\{s_i\}_{i=1,\ldots r-1}$ be the set
of simple reflections. Let $\widehat{\Si}_r= \Si_r \ltimes\Z^r$ be the
\textit{extended} affine symmetric group and let $s_0$ be the affine simple
reflection. Let $\Pi$ stand for the set of partitions and let $\Pi_r$ be the
set of partitions of length at most $r$. Elements of $\Pi^l$
for some $l \in \N$ will be called $l$-\textit{multipartitions}.
Finally, we will denote by
$\overline{Y}$ the Zariski closure of any subset $Y$ of an algebraic variety
$X$.
\section{Affine Hecke algebras and canonical bases} 
\paragraph{1.1} Consider the
Iwahori-Hecke algebra $\AH_r$ associated to $\widehat{\Si}_r$, i.e
the $\A$-algebra generated by elements $T_\sigma$, $\sigma \in \widehat{\Si}_r$
with relations
$$(T_{s_i}+1)(T_{s_i}-v^{-2})=0\qquad {for}\; i=0,\ldots,r-1,$$
$$T_{\sigma}T_{\gamma}=T_{\sigma\gamma}\qquad {if}\;
l(\sigma\gamma)=l(\sigma)l(\gamma).$$
We set $\tilde{T}_{\sigma}=v^{l(\sigma)}T_\sigma$ for every $\sigma\in 
\widehat{\Si}_r$.
\paragraph{}It is well-known that $\AH_r$ admits another presentation (the 
\textit{Bernstein} presentation) as the unital $\A$-algebra generated by
elements $T_i^{\pm 1}, X_j^{\pm 1}$ where $i \in [1,r-1]$, $j \in [1,r]$ with
the following relations
\begin{alignat*}{2}
&T_i\,T_i^{-1}=1=T_i^{-1}\,T_i,\qquad & \qquad &(T_i+1)(T_i-v^{-2})=0,\\
&T_i\,T_{i+1}\,T_i=T_{i+1}\,T_i\,T_{i+1},\qquad&\qquad&
|i-j|>1\Rightarrow T_i\,T_j=T_j\,T_i,\\
&X_i\,X_i^{-1}=1=X_i^{-1}\,X_i,\qquad&\qquad &X_i\,X_j=X_j\,X_i,\qquad\\
&T_i\,X_i\,T_i=v^{-2}X_{i+1},\qquad&\qquad 
&j\not= i,i+1\Rightarrow X_j\,T_i=T_i\,X_j.
\end{alignat*}
The isomorphism between the two presentations is such that
$T_{s_i}\mapsto T_i$ and $\tilde{T}_\lambda^{-1}\mapsto X_1^{\lambda_1}\cdots
X_r^{\lambda_r}$ if $\lambda=(\lambda_1,\ldots,\lambda_r)$ is
\textit{dominant}. The center of $\AH_r$ is $Z(\AH_r)=\A[X_1^{\pm 1},\ldots ,
X_r^{\pm 1}]^{\Si_r}$. Set 
$Z^-_r=\A[X_1^{-1},\ldots ,
X_r^{-1}]^{\Si_r}.$
\paragraph{1.2}For every $t,s \in \N$ define the left (resp. right)
representation of $\widehat{\Si}_t$ on $\Z^t$ of \textit{level} $s$ by
\begin{align*}
s_j \cdot (i_1, \ldots, i_t)&=(i_1, \ldots, i_{j+1}, i_j, \ldots i_t),
\qquad 1 \leq j <r,\\
\lambda \cdot (i_1, \ldots, i_t)&=(i_1+s\lambda_1,\ldots,i_t +s\lambda_t),
\qquad \lambda \in \Z^t
\end{align*}
and
\begin{align*}
(i_1, \ldots, i_t)\cdot s_j&=(i_1, \ldots, i_{j+1}, i_j, \ldots i_t),
\qquad 1 \leq j <r,\\
(i_1, \ldots, i_t)\cdot \lambda&=(i_1+s\lambda_1,\ldots,i_t +s\lambda_t),
\qquad \lambda \in \Z^t
\end{align*}
respectively. The set $\mathcal{A}^s_t=\{1 \leq i_1 \leq \cdots \leq i_t
\leq s\}$ is a fundamental domain for both actions. For each $\mathbf{i} \in
\mathcal{A}^s_t$ we set $\Si_\mathbf{i}=Stab\;\mathbf{i} \subset \Si_t$ and
denote by $\omega_{\mathbf{i}} \in \Si_\mathbf{i}$ the longest element. We also
let $\Si^\mathbf{i}$ be the set of all minimal length elements of the cosets
$\Si_\mathbf{i}\backslash \widehat{\Si}_t$.
\paragraph{1.3} Fix some $n \in \N^*$. For any $\mathbf{i},\mathbf{j} \in
\mathcal{A}^n_r$ and any
$\sigma \in \Si_{\mathbf{i}} \backslash
\widehat{\Si}_r/\Si_{\mathbf{j}}$ we set
$T_{\sigma}=\sum_{\delta \in \sigma} T_{\delta}$
and we let $\AH_{\mathbf{ij}} \subset \AH_r$ be the $\A$-linear span of the
elements $T_{\sigma}$ for $\sigma \in \Si_{\mathbf{i}}
\backslash \widehat{\Si}_r/\Si_{\mathbf{j}}$. Set $e_{\mathbf{i}}=\sum_{\delta
\in \Si_{\mathbf{i}}} T_{\delta}$. Then $\AH_{\mathbf{ij}}=e_{\mathbf{i}} \AH_r
e_{\mathbf{j}}$. Put
$$\widehat{\mathbf{S}}_{n,r}=\bigoplus_{\mathbf{i},\mathbf{j} \in
\mathcal{A}^n_r} \AH_{\mathbf{ij}}.$$
This space, equipped with the multiplication
$$e_{\mathbf{i}} h e_{\mathbf{j}} \bullet e_{\mathbf{k}} h' e_{\mathbf{l}}=
\delta_{jk} e_{\mathbf{i}} h e_{\mathbf{j}} h' e_{\mathbf{l}} \in
\AH_{\mathbf{il}}\qquad {for\;all}\; h,h' \in \AH_r$$
is called the \textit{affine q-Schur algebra}. It is proved in \cite{GV},
\cite{L3} that $\widehat{\mathbf{S}}_{n,r}$ is a quotient of the
\textit{modified} quantum affine algebra $\dot{\U}^-_v(
\widehat{\mathfrak{gl}}_n)$.

\paragraph{1.3}Set $\mathbf{T}_{n,r}=\bigoplus_{\mathbf{i} \in
\mathcal{A}^n_r} e_{\mathbf{i}}\AH_r$. For $\sigma \in
\Si_{\mathbf{i}}\backslash \widehat{\Si}_r$ we put $T_{\sigma}=
\sum_{\delta \in \sigma} T_\delta$. Then 
$\{\mathbf{T}_{\sigma}\},\;\sigma \in \Si_{\mathbf{i}}
\backslash \widehat{\Si}_r$ is an $\A$-basis of $e_\mathbf{i} \AH_r$. It will
be convenient to identify the element $\sigma$ with $\mathbf{i}
\cdot \sigma \in \Z^r$, so that $\{\mathbf{T}_p\},\;
p \in \Z^r$ is an $\A$-basis of $\mathbf{T}_{n,r}$.
\paragraph{}The algebra $\AH_r$ acts on $\mathbf{T}_{n,r}$ by multiplication
on the right, and $\widehat{\mathbf{S}}_{n,r}$ acts on $\mathbf{T}_{n,r}$ on
the left by
$$e_{\mathbf{i}} h e_{\mathbf{j}} \cdot e_{\mathbf{k}} h'=\delta_{\mathbf{jk}}
e_{\mathbf{i}}h e_{\mathbf{j}} h' \in e_{\mathbf{i}} \AH_r\qquad {for\;
every\;}h,h' \in \AH_r.$$
Let us denote these actions by $\rho_r:\widehat{\mathbf{S}}_{n,r} \to
\mathrm{End}\;(\mathbf{T}_{n,r})$ and $\sigma_r:\AH_r \to \mathrm{End}\;
(\mathbf{T}_{n,r})$. It is obvious that these two actions commute. The 
following result is a quantum and affine analogue of Schur-Weyl duality.
\begin{theo}[\cite{VV}] We have $\widehat{\mathbf{S}}_{n,r}=
\mathrm{End}_{\AH_r}(\mathbf{T}_{n,r})$. Moreover, we have $\AH_r=
\mathrm{End}_{\widehat{\mathbf{S}}_{n,r}}(\mathbf{T}_{n,r})$ if $n \geq r$.
\end{theo}

\paragraph{1.4}Let us now, following \cite{GV} and \cite{IM},
give the geometric realization
of the above Schur-Weyl duality. Let $\mathbb{L}=\F((z))$
and set $\mathbb{G}=GL_r(\mathbb{L})$. By definition, a \textit{lattice} in
$\mathbb{L}^r$ is a free $\F[[z]]$-submodule of rank $r$. Consider the variety
$X$ of sequences of lattices $(L_i)_{i \in \Z}$ such that
$$L_i \subset L_{i+1},\qquad \mathrm{dim}_{\F}(L_{i}/L_{i-1})=1,\qquad
L_{i+r}=z^{-1}L_i$$
(the \textit{affine flag variety of type $GL_r$}). Consider also the variety 
$Y$ of all $n$-step periodic flags in $\mathbb{L}^r$, i.e the set of all
sequences of lattices $(L_i)_{i \in \Z}$ such that
$$L_{i} \subset L_{i+1},\qquad L_{i+n}=z^{-1}L_i$$
(the \textit{affine partial flags variety}). The group $\mathbb{G}$ acts
(transitively) on $X$ and acts on $Y$ in obvious ways. Consider the diagonal
action of $\mathbb{G}$ on $X\times X$ and $Y \times Y$ respectively.
\paragraph{}It is well-known that the set of $\mathbb{G}$-orbits on
$X \times X$ is canonically identified with $\widehat{\Si}_r$. In order
to describe these $\mathbb{G}$-orbits we let $(e_1,\ldots,e_r)$
be a fixed $\mathbb{L}$-basis of $\mathbb{L}^r$
and set $e_{i+kr}=z^{-k}e_i$. Consider the right action of $\widehat{\Si}_r$
on $\Z^r$ of level $r$. To any element $\mathbf{x}$ in the orbit of
$\rho_r=(1,2,\ldots,r)$ we associate the flag $(L(\mathbf{x})_i)_{i \in \Z}$
defined by
$$L(\mathbf{x})_i=\prod_{p(j) \leq i} \F e_j,$$
where $p: \Z \to \Z$ is the bijection uniquely defined by $p(j)=\mathbf{x}_j$
if $1 \leq j \leq r$ and $p(j+r)=p(j)+r$. The $\mathbb{G}$-orbit
decomposition of $X \times X$ reads
$$X \times X=\bigsqcup_{\sigma \in \widehat{\Si}_r} X_\sigma$$
where $X_\sigma=\mathbb{G}\cdot (L(\rho_r\cdot\sigma),L(\rho_r))$. Similarly,
to each $\mathbf{i} \in \Z^r$
we associate the map $p: \Z \to \Z$ uniquely defined by $p(j)=\mathbf{i}_j$
if $1 \leq j \leq r$ and $p(j+r)=p(j)+n$. Consider the flag
$$L(\mathbf{i})_i=\prod_{\mathbf{i}(j) \leq i} \F e_j.$$
Then $Y=\bigsqcup_{\mathbf{i} \in \mathcal{A}^n_{r}} Y_{\mathbf{i}}$ where
$Y_{\mathbf{i}}=\mathbb{G}\cdot (L(\mathbf{i}))$ and
$$Y_{\mathbf{i}} \times Y_\mathbf{j}=\bigsqcup_{\sigma \in
\Si_{\mathbf{i}}\backslash \widehat{\Si}_r/\Si_{\mathbf{j}}}
Y_{\sigma}$$
where $Y_{\sigma}
=\mathbb{G}\cdot(L(\mathbf{i}\cdot \sigma),L(\mathbf{j}))$ and where the right
action of $\widehat{\Si}_r$ on $\Z^r$ is now of level $n$.

\paragraph{}Let $\C_{\mathbb{G}}(X \times X)$ (resp. $\C_{\mathbb{G}}
(Y \times Y)$) be the space of complex-valued $\mathbb{G}$-invariant functions
on $X \times X$ (resp. on $Y \times Y$) which are supported on finitely many 
orbits. The convolution product endows these
spaces with an associative algebra structure. We let
$\mathbf{1}_{\mathcal{O}} \in
\C_{\mathbb{G}}(X \times X)$ (resp. $ \mathbf{1}_{\mathcal{O}}\in 
\C_{\mathbb{G}}
(Y \times Y)$) be the
characteristic function of a $\mathbb{G}$-orbit $\mathcal{O} \subset X \times
X$(resp. $\mathcal{O}\subset Y \times Y$).
\begin{theo}[\cite{IM},\cite{VV}]\hfill
\begin{enumerate}
\item[i)]The linear map $(\AH_r)_{|v=q^{-1}} \to \C_{\mathbb{G}}
(X \times X)$ defined by $T_\sigma \mapsto \mathbf{1}_{X_\sigma}$ is an 
algebra isomorphism.
\item[ii)] The linear map $(\widehat{S}_{n,r})_{|v=q^{-1}}\to \C_{\mathbb{G}}
(Y \times Y)$ such that $T_{\sigma}\mapsto
\mathbf{1}_{Y_{\sigma}}$
is an algebra isomorphism.
\end{enumerate}
\end{theo}
\paragraph{}Now consider the diagonal action of $\mathbb{G}$ on $Y \times X$.
The collection of orbits are parametrized by $\Z^r$: to
$\mathbf{i} \in \Z^r$ corresponds the orbit $\mathcal{O}_\mathbf{i}$
of the pair
$(L(\mathbf{i}),L(\rho_r))$. The algebras $\C_{\mathbb{G}}(X \times X)$ and
$\C_{\mathbb{G}}(Y \times Y)$ act by convolution on 
$\C_{\mathbb{G}}(Y \times X)$ on the right and on the left respectively.
\begin{theo}[\cite{VV}] The map $(\mathbf{T}_{n,r})_{|v=q^{-1}}\to
\C_{\mathbb{G}}(Y \times X)$ such that $e_\mathbf{i} \mapsto
\mathbf{1}_{\mathcal{O}_\mathbf{i}}$ for $\mathbf{i} \in \mathcal{A}^n_r$
extends uniquely to an isomorphism of
$(\widehat{\mathbf{S}}_{n,r})_{|v=q^{-1}}\times (\AH_r)_{|v=q^{-1}}$-modules.
\end{theo}

\paragraph{1.5}Let $u \mapsto \overline{u}$ be the semilinear involution of
$\AH_r$ defined by $\overline{T}_\sigma=T^{-1}_{\sigma^{-1}}$ for all $\sigma$.
For each $\sigma \in \widehat{\Si}_r$ there exists a unique element $\c_\sigma
\in \AH_r$ such that
$$\mathrm{i)} \;\;\overline{\c_\sigma}=\c_{\sigma},\qquad \mathrm{ii)}\;\;
\c_\sigma=\tilde{T}_\sigma
+\sum_{\delta<\sigma}c_{\delta,\sigma}(v)\tilde{T}_\delta,\;\quad c_{\delta,
\sigma}(v)\in v\S.$$
The polynomial $c_{\sigma,\delta}(v)$ is the affine Kazhdan-Lusztig polynomial
of type $\tilde{A}_{r-1}$ associated to $\sigma$ and $\delta$ (this polynomial
is denoted by $h_{\sigma,\delta}$ in Soergel's notation \cite{Soergel}).
\paragraph{}For $\sigma \in \widehat{\Si}_r$ and $L \in X$ let $X_{\sigma,L}$
be the fiber of the first projection $X_{\sigma} \to X$. 
Then $X_{\sigma,L}$ is the set of $\F$-points of an algebraic
variety of dimension $l(\sigma)$ whose isomorphism class is independent of $L$.
Then
$$\mathbf{c}_\sigma=\sum_{i,\delta} v^{-i+l(\sigma)-l(\delta)}\mathrm{dim}\;
\mathcal{H}^i_{X_{\delta,L}}(IC_{X_{\sigma,L}}) \tilde{T}_\delta$$
where $IC_{X_{\sigma,L}}$ denotes the intersection cohomology complex
associated to $X_{\sigma,L}$ and where $\mathcal{H}^i$ stands for
local cohomology.
\paragraph{}Similarly, let $\mathbf{i},\mathbf{j} \in \mathcal{A}^n_r$
and let
$\sigma \in
\Si_{\mathbf{i}}\backslash \widehat{\Si}_r/\Si_{\mathbf{j}}$. 
Denote by $Y_{\sigma,\mathbf{i}}$ the
fiber above $(L(\mathbf{i}))$ of the projection of $Y_{\sigma} \to
Y$ on the first component. This is the set of $\F$-points of an algebraic
variety of dimension, say $y(\sigma)$ (an explicit formula for
$y(\sigma)$ can be found in \cite{L3}). Put
$\tilde{T}_{\sigma}=v^{y(\sigma)}T_{\sigma}$.
For every $\sigma \in \Si_{\mathbf{i}}\backslash
\widehat{\Si}_r/\Si_{\mathbf{j}}$ set
$$\mathbf{c}_{\sigma}=\sum_{i,\delta}
v^{-i+y(\sigma)-
y(\delta)} \;\mathrm{dim}\;\mathcal{H}^i_{
Y_{\delta,\mathbf{i}}}
(IC_{Y_{\sigma,\mathbf{i}}}) \tilde{T}_{\delta}.$$
It is clear that
$\overline{\AH_{\mathbf{ij}}}=\AH_{\mathbf{ij}}$. Define a semilinear
involution $\tau : \AH_{\mathbf{ij}} \to \AH_{\mathbf{ij}}$ by $\tau (u)=
v^{-2l(\omega_\mathbf{j})} \overline{u}$.
The elements $\{\mathbf{c}_{\sigma}\}$ for all $\mathbf{i},\mathbf{j} \in
\mathcal{A}^n_r$ form the
canonical basis of $\widehat{\mathbf{S}}_{n,r}$ and are characterized by the
following two properties :
$$\mathrm{i)}\;\tau(\mathbf{c}_{\sigma})=
\mathbf{c}_{\sigma},\qquad
\mathrm{ii)}\;\mathbf{c}_{\sigma}=\tilde{T}_{\sigma}+
\sum_{\delta
<\sigma} c_{\delta,\sigma}(v)
\tilde{T}_{\delta},\qquad
c_{\delta,\sigma}(v) \in v\S.$$
\paragraph{1.6} Let $s,t \in \N^*$. For $\mathbf{i} \in \mathcal{A}^s_t$ and 
$x \in \mathbf{i} \cdot \widehat{\Si}_t$ set $\langle x|=e_\mathbf{i}
\tilde{T}_a$ where $\mathbf{i}\cdot a=x$ and $a \in \Si^\mathbf{i}$. The set
$\{\langle x|, x \in \mathbf{i} \cdot \widehat{\Si}_t\}$ is an $\A$-basis
of the space $e_\mathbf{i} \AH_t$. Define a semilinear involution $u \mapsto
\overline{u}$ of $e_\mathbf{i} \AH_t$ by $\overline{e_\mathbf{i}x}=
e_{\mathbf{i}}\overline{x}$. There exists a unique $\A$-basis
$\{\mathbf{c}^-_{x},
x \in \mathbf{i} \cdot \widehat{\Si}_t\}$ of $e_\mathbf{i} \AH_t$ such that
$$\mathrm{i)}\;\overline{\mathbf{c}^-_{x}}=
\mathbf{c}^-_x,\qquad
\mathrm{ii)}\;\mathbf{c}^-_{x}=\langle x|+
\sum_{y} P^-_{y,x}
\langle y|,\qquad
P^-_{y,x}\in v^{-1}\Z[v^{-1}].$$
The polynomials $P^-_{y,x}$ are \textit{parabolic affine
Kazhdan-Lusztig polynomials} introduced by Deodhar \cite{Deo}. These
polynomials are (up to a sign) denoted by $\overline{n}_{a_y,a_x}$ in
Soergel's notation, where $a_x,a_y \in \Si^\mathbf{i}$ are such that
$x=\mathbf{i}\cdot a_x$, $y=\mathbf{i}\cdot a_y$. 
\section{The main result}
\paragraph{2.1} Let $\Gamma$ be Macdonald's ring of symmetric polynomial
in the variables $y_i$, $i \in \Z$, defined over $\A$ (see \cite{Mac}). Let
$\Gamma_r=\A[y_1,\ldots,y_r]^{\Si_r}$. Let $s_\lambda \in \Gamma_r$ be the
Schur polynomial associated to $\lambda \in \Pi_r$.
\paragraph{}Fix some $n \in \N$ and let $\mathbf{i} \in
\mathcal{A}^n_r$. From $s_\lambda(X_1^{-1},\ldots,X_r^{-1}) \in Z(\AH_r)$ it
follows that
$e_\mathbf{i} s_\lambda(X_1^{-1},\ldots,X_r^{-1}) \in \AH_{\mathbf{ii}}$.
Define polynomials $J^\mathbf{i}_{\lambda,\sigma} \in \Z[v,v^{-1}]$
by the relation
$$e_\mathbf{i} s_\lambda(X_1^{-1},\ldots,X_r^{-1})=(-v)^{(n-1)|\lambda|}
\sum_{\sigma \in
\Si_{\mathbf{i}} \backslash \widehat{\Si}_r / \Si_{\mathbf{i}}}
J^\mathbf{i}_{\lambda,\sigma} \mathbf{c}_{\sigma}.$$
In this section we give an explicit expression for
$J^\mathbf{i}_{\lambda,
\sigma}$ involving (parabolic) affine Kazhdan-Lusztig polynomials of
type $A$. 
\paragraph{Remark.}It is clear that (up to a power of $v$)
$J^\mathbf{i}_{\lambda,\sigma}$ depends only on
$\Si_\mathbf{i}$ rather than on $\mathbf{i}$. In particular, any parabolic
subgroup $\Si_{i_1}\times \cdots \times \Si_{i_t}$ occurs as
$\Si_\mathbf{i}$ for some $\mathbf{i} \in \mathcal{A}^n_r$ as soon as $n \geq
t$.

\paragraph{2.2} We first make some preliminary definitions.
 We will represent a partition $\lambda$ by its
associated Young diagram in the usual fashion. We will consider diagrams
where the $(i,j)$-box has content $i-j+r_0\;mod\;n$ for some fixed $r_0 \in
\Z/n\Z$ and call the resulting tableau the \textit{partition $\lambda$ with
residue $r_0$}. We will say that a box with content $j \in \Z/n\Z$ can be
\textit{added} to the partition $\lambda$ with residue $r_0$ if there exists a
partition $\lambda'$ with residue $r_0$ such that $\lambda'/\lambda$
is a single box with content $j$. For example, when $n=3$, the partition
$\lambda=(421)$ with residue $1$ is
$$
\centerline{\epsfbox{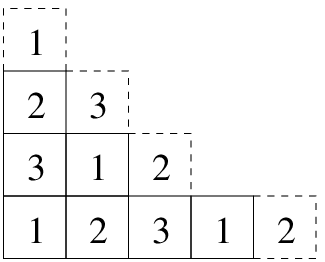}}
$$
and the dotted lines correspond to addable boxes.
\paragraph{} To each $\mathbf{p} \in (\Z^+)^r$ and $\mathbf{i} \in
\mathcal{A}^n_r$ we
associate a multipartition (with residues)
$\mathcal{M}_\mathbf{i}(\mathbf{p})$. 
First, we attach a diagram (not a partition!)
$$D_{\mathbf{p}}=\{(i,j)\;|\;0<j\leq \mathbf{p}_i\}\subset \Z/r\Z 
\times \Z^+,$$
where we fill the $(i,j)$-box with the content $\mathbf{i}_i+
\mathbf{p}_i-j\;mod\;n$.
\noindent
\paragraph{}\textit{Example 1.} Suppose $r=n=5$, $\mathbf{i}=(1,2,3,4,5)$ and
$\mathbf{p}=(4,3,4,3,5)$.
Then $D_{\mathbf{p}}$ is
$$
\centerline{\epsfbox{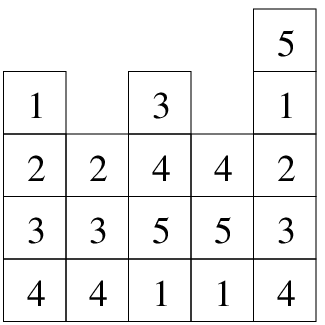}}
$$
Now consider the horizontal slices $s_k=D_{\mathbf{p}}
\cap (\Z/r\Z \times \{k\})$ and let
$k_0$ be maximal such that $s_{k_0}\neq \emptyset$. We construct the
multipartition with residues $\mathcal{M}_\mathbf{i}(\mathbf{p})$
by successively
adding the boxes from $s_{k_0},\ldots,s_1$ in the following way. Set
$\mathcal{M}^{k_0+1}=\emptyset$. Suppose $\mathcal{M}^i=(\lambda^{(1)}_i,\ldots
,\lambda^{(t)}_i)$ is known. Then $\mathcal{M}^{i-1}=(\lambda^{(1)}_{i-1},
\ldots,\lambda^{(r)}_{i-1})$ is obtained from
$\mathcal{M}^i$ by adding the boxes from $s_i$ (possibly creating new
partitions) in such a way that
\begin{enumerate}
\item[i)] For every $1 \leq v \leq r$, $\lambda^{(v)}_{i-1} / \lambda^{(v)}_i$
is a skew tableau with \textit{at most} one box in each row,
\item[ii)] $\mathcal{M}^{i-1}$ is maximal for the following order :\\
$(\lambda_{i-1}^{(1)},\lambda_{i-1}^{(2)},\ldots) \geq (\mu_{i-1}^{(1)},
\mu_{i-1}^{(2)},\ldots)\;if\;there\;exists\;w\;such\;that$
$$\lambda_{i-1}^{(l)}=\mu_{i-1}^{(l)}\;\;{for\;} 1\leq l <w\qquad
{and}\qquad\lambda_{i-1}^{(w)} \geq \mu_{i-1}^{(w)},$$
where $\geq$ stands for the usual dominance order of partitions,
\item[iii)] If several new partitions appear in $\mathcal{M}^{i-1}$ then they
are in increasing order of their residue.
\end{enumerate}
Set $\mathcal{M}_\mathbf{i}(\mathbf{p})=\mathcal{M}^1$.
We note that condition iii) above is
 not essential for the rest of the paper and here only to fix notations.
\noindent
\paragraph{}\textit{Examples.} i) Let $r=n$ and $\mathbf{i}=(1,\ldots,r)$.
Suppose that $\mathbf{p}$
is \textit{antidominant} up to cyclic permutation, i.e there exists 
$i \in \Z/r\Z$ such that $\mathbf{p}_i \geq \mathbf{p}_{i-1} \geq
\cdots \geq \mathbf{p}_{i+1}$.
Let $\lambda$ be the associated partition. Then
$\mathcal{M}_\mathbf{i}(\mathbf{p})$ consists of the single partition
$\lambda$ with residue~$i$.\\
ii) Consider $r,n$, $\mathbf{i}$ and $\mathbf{p}$ as in example 1. Then the
algorithm for computing $\mathcal{M}_\mathbf{i}(\mathbf{p})$ runs as follows :
$$
\centerline{\epsfbox{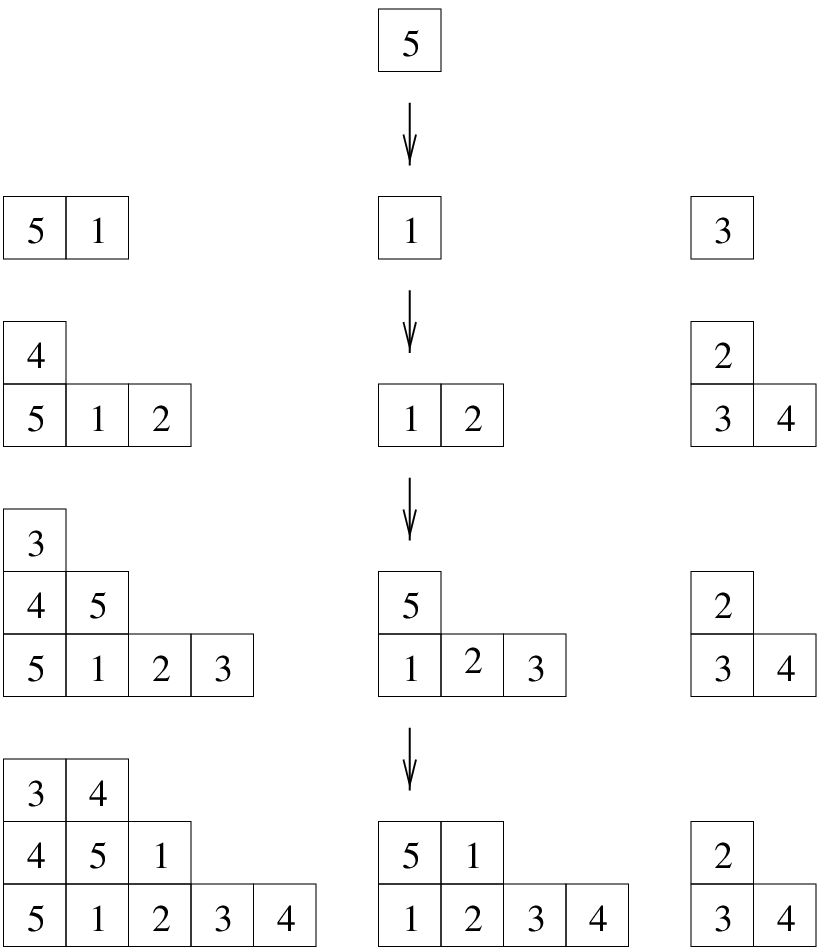}}
$$
 For $\sigma \in \Si_{\mathbf{i}} \backslash
\widehat{\Si}_r / \Si_{\mathbf{i}}$ we let $\sigma_0$ be the
longest element in $\sigma$ and we set
$\mathbf{i}\bullet \sigma =\mathbf{i}\cdot
\sigma_0-\mathbf{i} \in \Z^r$. For each $\sigma$ such that
$\mathbf{i}\bullet\sigma \in (\Z^+)^r$ we set $\mathcal{M}(\sigma)=
\mathcal{M}_\mathbf{i}(\mathbf{i}\bullet\sigma)$. Write
$\mathcal{M}(\sigma)=(\sigma^{(1)},\ldots,\sigma^{(l)})$ where $\sigma^{(l)}
\neq \emptyset$, and $\mathbf{r}_\sigma=(r_1,\ldots,r_l)$ where $r_i \in 
\Z/n\Z$ is the residue of $\sigma^{(i)}$.

\paragraph{2.3} Let $l \in \N$. Let $(\sigma^{(1)},\ldots,\sigma^{(l)})$,
 $(\mu^{(1)},\ldots,\mu^{(l)})$
be any $l$-multipartitions and let $\mathbf{r}=(r_1,\ldots,r_l) \in
(\Z/n\Z)^l$. Choose some $\mathbf{s}=(s_1,\ldots,s_l) \in \Z^l$ such that
$s_i \equiv r_i \;(\mathrm{mod}\;n)$.
For $i=1,\ldots l$ and $j \in \N$ we set $u_j^{(i)}=s_i+\sigma^{(i)}_j+1-j$ and
$v_j^{(i)}=s_i+\mu^{(i)}_j+1-j$. Consider, for $t \gg 0$
$$\mathbf{u}=(u^{(1)}_1,\ldots,u^{(1)}_t,u^{(2)}_1,\ldots, u^{(2)}_t,\ldots,
u^{(l)}_t),$$
$$\mathbf{v}=(v^{(1)}_1,\ldots,v^{(1)}_t,v^{(2)}_1,\ldots, v^{(2)}_t,\ldots,
v^{(l)}_t).$$
Finally, we put
$$\mathbf{P}^{-,\mathbf{s}}_{(\mu^{(1)},\ldots,\mu^{(l)}),
(\sigma^{(1)},\ldots,\sigma^{(l)})}=P^-_{\mathbf{v},\mathbf{u}}.$$
Now let $\mathbf{s}$ be in the asymptotic range $s_1 \gg s_2 \gg \cdots \gg 
s_l$ and set 
$$\mathbf{P}^{-,\mathbf{r}}_{(\mu^{(1)},\ldots,\mu^{(l)}),
(\sigma^{(1)},\ldots,\sigma^{(l)})}=P^-_{\mathbf{v},\mathbf{u}}.$$
This polynomial is independent
of the choices of $\mathbf{s}$ and $t$ in the given asymptotic range (this
follows for instance from \cite{U} Section 4 and \cite{S}, Theorem 4.1).
These can be thought of as some ``stabilization'' of polynomials 
$P^-_{\mu+\rho_r,\sigma+\rho_r}$ of type $\tilde{A}_{r}$ as
$r$ tends to infinity (see \cite{LT}). Moreover, it is easy to see that when
$l=1$, $\mathbf{P}^{-,\mathbf{r}}$ is independent of $\mathbf{r}$ and we will
omit it.

\paragraph{2.4}For any multipartition $\mu=(\mu^{(1)},\ldots,\mu^{(l)})$
and $\mathbf{r}=(r_1,\ldots,r_l) \in (\Z/n\Z)^l$ we set $\mu'=((\mu^{(l)})'
,\ldots,(\mu^{(1)})')$ and $\mathbf{r}'=(-r_l,\ldots,-r_1)$.
\paragraph{2.5}The following is the main result of this paper, and will be
proved in Section~5.
\begin{theo} We have 
$$e_{\mathbf{i}}s_\lambda(X_1^{-1},\ldots,X_r^{-1})=(-v)^{(n-1)|\lambda|}
\sum_{\sigma\;|\mathbf{i}\bullet\sigma \in (\Z^+)^r} 
J^\mathbf{i}_{\lambda,\sigma} \mathbf{c}_\sigma$$
where
\begin{equation*}\label{E:Main}
J^\mathbf{i}_{\lambda,\sigma}=\sum_{\stackrel{\nu_1,\ldots,\nu_l}{\mu_1,
\ldots,\mu_l}} c^\lambda_{\mu_1,
\ldots,\mu_l} v^{\sum (b-1)|\mu_b|}
\mathbf{P}^-_{\nu_1,n\mu_1'}\cdots \mathbf{P}^-_{\nu_l,n\mu_l'}
{\mathbf{P}^{-,\mathbf{r}_\sigma'}_{\nu,
\mathcal{M}(\sigma)'}}
\end{equation*}
and $\nu=(\nu_1,\ldots,\nu_l)$. Here
$c^\lambda_{\mu_1,\ldots,\mu_l}$ is the (generalized) 
Littlewood-Richardson coefficient.
\end{theo}
\paragraph{}\textit{Examples.} i) Suppose that $n=1$. Then $\mathbf{i}=(1^r)$
and $\Si_\mathbf{i}=\Si_r$. Moreover, $\Si_r \backslash \widehat{\Si}_r / 
\Si_r =\Pi_r$ and for $\sigma \in \Pi_r$ we have $\mathcal{M}(\sigma)=\sigma$
and $l=1$.
Hence the above theorem reduces to $J^\mathbf{i}_{\lambda,\sigma}=
\sum_\nu \mathbf{P}^-_{\nu,\lambda'}\mathbf{P}^-_{\nu,\sigma'}=
\delta_{\lambda,\sigma}$, i.e
$$\big(\sum_{w \in \Si_r}T_w\big)s_\lambda(X_1^{-1},\ldots,X_r^{-1})=
\mathbf{c}_\lambda,$$
 in accordance with \cite{L1}.\\
ii) Let $r=n$ and $\mathbf{i}=\rho$ (i.e $\Si_\mathbf{i}=\{1\}$). Let
$\lambda=(1^l)$, $l \leq r$ be a minuscule weight. Then in the above expression
for $J^\mathbf{i}_{\lambda,\sigma}$ the only nonzero terms correspond to the
case when $\mu_i$ is also minuscule for all $i$. We obtain an expression
for $s_\lambda(X_1^{-1},\ldots,X_r^{-1})$ analogous to Theorem 1.1
in \cite{Hai2} for $G=GL(r)$ (but which involves Kazhdan-Lusztig polynomials
rather than $R$-polynomials). Note that \cite{Hai1}, Proposition 5 also easily
follows from the above theorem.

\section{Hall algebra of a cyclic quiver}
\paragraph{3.0 Notations.} In this section we fix a positive integer $n$. Let
$(\epsilon_i)$, $i \in \Z/n\Z$ be the canonical basis of
$\N^{\Z/n\Z}$. For $i \in \Z/n\Z$ and $l \in \N^*$, define the \textit{cyclic
segment} $[i;l)$ to be the image of the projection to $\Z/n\Z$ of the segment
$[i_0,i_0+l-1] \subset \Z$ for any $i_0 \equiv i\;(\mathrm{mod}\;n)$. A
{\textit{cyclic multisegment} is a linear combination $\m=\sum_{i,l} a^l_i
[i;l)$ of cyclic segments with coefficients $a^l_i \in \N$. Let $\mathcal{M}$
be the set of cyclic multisegments. For $\m \in \mathcal{M}$ we set
$\mathrm{dim}\;\m=\sum a^l_i (\epsilon_i + \cdots + \epsilon_{i+l-1})$.
Note that $\mathcal{M}$ is canonically
isomorphic to $\Pi^n$ : to $\m=\sum a_i^l [i;l)$ we associate the 
multipartition $(\lambda^{(1)},\ldots,\lambda^{(n)})$ with $\lambda^{(i)}=
(1^{a_i^1}2^{a_i^2}\ldots)$. 

\paragraph{3.1} Let $Q$ be the quiver of type $\tilde{A}_{n-1}$, i.e the
oriented graph with vertex set $I=\Z/n\Z$ and edge set $\Omega=\{({i},{i+1}), 
i \in I\}$. For any $I$-graded $\mathbb{F}$-vector space $V=\bigoplus_{i \in 
I} V_i$, let $E_V \subset \bigoplus_{(i,j) \in \Omega} \mathrm{Hom}\;(V_i,
V_j)$ denote the space of nilpotent representations of $Q$.
The group $G_V=\prod_{i\in I} 
GL(V_i)$ acts on $E_V$ by conjugation. For each $i \in I$ there exists a 
unique simple $Q$-module $S_i$ of dimension $\epsilon_i$, and for each pair 
$(i,l) \in I \times \N^*$ there exists a unique (up to isomorphism) 
indecomposable $Q$-module $S_{i;l}$ of length $l$ and tail $S_i$. 
Furthermore, every nilpotent $Q$-module $M$ admits an essentially unique 
decomposition
\begin{equation}\label{E:1.1}
M \simeq \bigoplus_{i,l} a_i^lS_{i;l}.
\end{equation}
We denote by $\overline{\mathbf{m}}$ the isomorphism class of $Q$-modules 
corresponding (by (\ref{E:1.1})) to the multisegment $\mathbf{m}=\sum_{i,l} 
a_i^l [i;l)$. For $\m \in \mathcal{M}$ with
$\mathrm{dim}\;\m=\mathbf{d}$ and $V_\mathbf{d}$ an $I$-graded vector space
of dimension $\mathbf{d}$, we let $O_\m \subset E_{V_\mathbf{d}}$ be the $G_{V_\mathbf{d}}$-orbit consisting of representations in the class $\overline{\mathbf{m}}$, and we let $\mathbf{1}_\mathbf{m} \in \C_G(V_{\mathbf{d}})$ be the characteristic function of $O_\m$. Finally, we set $\f_\m=q^{-\mathrm{dim}\;O_\m} \mathbf{1}_\m$. We will write $\m < \n$ if $\mathcal{O}_\m \subset
\overline{\mathcal{O}_\n}$.
\paragraph{3.2} Set $\U^-_n=\bigoplus_{\mathbf{d}} \C_G(E_{V_\mathbf{d}})$. Note that, by definition, $(\f_\m)_{\m \in \mathcal{M}}$ is a $\C$-basis of $\U^-_{n}$. The space $\U^-_n$ is endowed with the structure of a (Hall) algebra
(see \cite{L1}). We use the definitions of \cite{VV}, \cite{S}. Moreover, the
structure constants for this algebra are polynomials in $q$, and one can consider $\U^-_n$ as an $\A$-algebra with $q=v^{-1}$. The algebra $\U^-_n$ is naturally $\N^{\Z/n\Z}$-graded and we denote by $\U^-_n[\mathbf{d}]$ the component of degree $\mathbf{d}$.
Let $\U_v(\widehat{\mathfrak{sl}}_n)$
denote the Lusztig integral form of the quantum affine algebra of type $\tilde{A}_{n-1}$ and let $e_i^{(l)},k_i,f_i^{(l)}$, $i\in I$, $l \in \N$ be the
divided powers of the standard Chevalley generators. Let $\uasn$ be the subalgebra of $\U_v(\widehat{\mathfrak{sl}}_n)$ generated by $f^{(l)}_i, i\in I$, $l \in \N^*$. It is known that the map $f^{(l)}_i \mapsto
\mathbf{f}_{l\epsilon_i}$ extends to an embedding of the algebras $\uasn 
\hookrightarrow \U^-_n$.
\paragraph{3.3} For $\mathbf{m}\in \mathcal{M}$, set
\begin{equation}\label{E:canHall}
\mathbf{b}_{\mathbf{m}}=\sum_{i,\mathbf{n}}v^{-i+\mathrm{dim}\,O_{\mathbf{m}}-\mathrm{dim}\,O_{\mathbf{n}}}\mathrm{dim}\,\mathcal{H}^i_{O_{\mathbf{n}}}(IC_{O_{\mathbf{m}}}) \mathbf{f}_{\mathbf{n}},
\end{equation}
where $\mathcal{H}^i_{O_{\mathbf{n}}}(IC_{O_{\mathbf{m}}})$ is the
stalk over a point of $O_{\mathbf{n}}$ of the ith intersection
cohomology sheaf of the closure $\overline{O}_{\mathbf{m}}$ of
$O_{\mathbf{m}}$. Then $\mathbf{B}=\{\mathbf{b}_{\mathbf{m}}\}$ is the
canonical basis of $\U^-_n$, introduced in \cite{VV}.

\paragraph{3.4}Let $L,L' \in Y$ be two $n$-step periodic flags in
$\mathbb{L}^r$ satisfying
$L' \subset L$. Following Lusztig (see \cite{L2},\cite{GV}) we associate to
such a pair
a nilpotent representation of $\tilde{A}_{n-1}$ of graded dimension
$(\mathrm{dim}_\F (L_i/L'_i))_{\bar{i} \in \Z/n\Z}$. Let us denote by
$L/L'$ this $\tilde{A}_{n-1}$-module. Set
$$a(L',L)=\sum_{i=1}^n \mathrm{dim}_\F(L_i/L'_i)(\mathrm{dim}_\F(L'_{i+1}/L'_i)
-\mathrm{dim}_\F(L_i/L'_i)).$$
Define a map $\Theta:\;\U^-_n \to \widehat{\mathbf{S}}_{n,r}$ by
$$\Theta(f)(L',L)=q^{-a(L',L)}f(L/L')\quad\mathrm{if}\;L'\subseteq L$$
and $\Theta(f)(L,L')=0$ if $L' \not\subseteq L$.
\paragraph{}In order to describe $\Theta$, we consider the following
parametrization of the collection of $\mathbb{G}$-orbits in $Y \times Y$. 
Let $M_{r,n}$ be the set of $\Z \times \Z$-matrices $\mathbf{s}=
(s_{ij})_{i,j \in \Z}$ with entries in $\N$ such that $s_{i+n,j+n}=s_{i,j}$
and $\sum_{j}\sum_{i=1}^n s_{ij}=r.$
To each such $\mathbf{s} \in M_{r,n}$ we associate the $\mathbb{G}$-orbit
$Y_{\mathbf{s}}$ whose elements are the pairs $(L,L')$ for which
$$s_{ij}=\mathrm{dim}_\F\; \bigg( \frac{L_i \cap L'_j}{(L_i \cap L'_{j-1}) +
(L_{i-1} \cap L_j')}\bigg).$$
For $\mathbf{i},\mathbf{j}\in \mathcal{A}^n_{r}$ we denote by $M_{\mathbf{i}
\mathbf{j}}$ the set of all $\mathbf{s}$ such that $Y_{\mathbf{s}}\subset
Y_{\mathbf{i}} \times {Y}_{\mathbf{j}}$. It is easy to see that
$$M_{ij}=\{\mathbf{s}\in M_{r,n}\;|\; \sum_j s_{ij}=\#\mathbf{i}^{-1}(i),\;
\sum_{i}s_{ij}=\#\mathbf{j}^{-1}(j)\}.$$
In particular, $M_{\mathbf{ij}}$ is naturally identified with
$\Si_{\mathbf{i}}\backslash \widehat{\Si}_r/\Si_{\mathbf{j}}$.
\paragraph{} Let us associate to each $\m=\sum a_i^l [i;l)$ the
matrix $(m_{i,j})\in \bigcup_r M_{r,n}$ with $m_{i,j}=a_i^{j-i+1}$.
The set 
$$M^+=\{(m_{i,j})_{i,j \in \Z}\:|\; m_{i+n,j+n}=m_{i,j},\;i>j \Rightarrow
m_{i,j}=0\}$$
is then identified with $\mathcal{M}$. If $\mathbf{i} \in
\mathcal{A}^n_{r}$ and $\m \in M^+$ we let $\m^\mathbf{i} \in 
\bigcup_\mathbf{j}M_{\mathbf{ij}}$ be the matrix whose $(i,j)$th
entry is
$$\delta_{ij}(\#\mathbf{i}^{-1}(j+1)-\sum_{k \leq j}m_{kj})+(1-\delta_{ij})
m_{i+1,j}.$$
\begin{prop}[\cite{VV}] The map $\Theta:\;\U^-_n \to
\widehat{\mathbf{S}}_{n,r}$ is an algebra morphism satisfying 
$\Theta(\overline{u})=\tau(\Theta(u))$ for every $u \in \U^-_n$. Furthermore,
$$\Theta(\mathbf{f}_\mathbf{m})=\sum_{\mathbf{i}\;|\;\mathbf{m}^\mathbf{i}
\in M^+} \tilde{T}_{\mathbf{m}^\mathbf{i}},\qquad 
\Theta(\mathbf{b}_\mathbf{m})=\sum_{\mathbf{i}\;|\;\mathbf{m}^\mathbf{i}
\in M^+} \mathbf{c}_{\mathbf{m}^\mathbf{i}}.$$
\end{prop}
\paragraph{}It follows from the above Proposition that $\mathbf{T}_{n,r}$
is endowed with a canonical $\U^-_n$-module structure.
\paragraph{3.5}Let $e'_i$, $i \in \Z/n\Z$ be the adjoint of the left
multiplication by $\mathbf{f}_i$. Set $\mathbf{R}=\bigcap_{i} \mathrm{Ker}\;
e'_i\subset \U^-_n$. Let us identify the ring of symmetric polynomials
$\Gamma_r$ with $Z^-_r$ by $y_i \mapsto X_i^{-1}$.
\begin{theo}[\cite{S}] The vector space $\mathbf{R}$ is a graded central
subalgebra of $\U^-_n$ and the multiplication map induces an isomorphism
$\uasn \otimes_\A \mathbf{R} \stackrel{\sim}{\to} \gaga$. Moreover there
exists surjective algebra morphisms $i_r:\mathbf{R} \to Z^-_r$ and an algebra
isomorphism $i:\mathbf{R} \to \Gamma$ such that
$$\rho_r \circ \Theta = \sigma_r \circ i_r, \qquad\;i=
\underset{\longleftarrow}{\mathrm{lim}}\; i_r.$$
\end{theo}
\paragraph{}Let $s_\lambda \in \Gamma$ be the Schur polynomial associated to
$\lambda \in \Pi$, and set $a_\lambda=i^{-1}(s_\lambda)$. Then $i_r(a_\lambda)
=s_\lambda(X_1^{-1},\ldots,X_r^{-1})$ for any $r \geq l(\lambda)$. For $\m \in \mathcal{M}$, define
polynomials $J^\lambda_\m \in \Z[v,v^{-1}]$ by
\begin{equation}\label{E:0def}
a_\lambda=\sum_\m J_{\lambda,\m}\b_\m.
\end{equation}
\begin{cor} For any $r \in \N$ and $\mathbf{i} \in \mathcal{A}^n_{r}$ we have
\begin{equation}\label{E:cor2.9}
\mathbf{e}_\mathbf{i} s_\lambda(X_1^{-1},\ldots X_r^{-1})
=\sum_{\mathbf{m}\;|\mathbf{m}^{\mathbf{i}}\in M^+} J_{\lambda,\mathbf{m}}
\mathbf{c}_{\mathbf{m}^\mathbf{i}}.
\end{equation}
\end{cor}
\noindent
\textit{Proof.} This follows by applying $a_\lambda$ to $\mathbf{e}_\mathbf{i}
\cdot 1 \in \mathbf{T}_{n,r}$, and using Theorem 3.5 and Proposition 3.4.\qed

\paragraph{Remarks.} i) Let us consider the case $n=1$ and $\mathbf{i}=(1^r)$. Then $\mathcal{M}=\Pi$
and $\U^-_1=\mathbf{R}\stackrel{i}{\simeq}\Gamma$, and it is known that $i$ 
identifies the Poincar\'e-Birkhoff-Witt basis element $\mathbf{f}_\lambda$
with the Hall-Littlewood polynomial $P_\lambda$ (see \cite{Mac}, Chap. III). In
particular, $K^\lambda_\mu(v)$ is the Kostka-Foulkes polynomial and from
(\ref{E:cor2.9}) we recover the well-known result of Lusztig
(\cite{L1}) concerning the Satake isomorphism
$$(\sum_{\sigma \in \Si_r} T_\sigma)s_\lambda(X_1^{-1},\ldots,X_r^{-1})
=\sum_{\mu \in \Pi} K^\lambda_\mu(v) \tilde{T}_{\Si_r \mu \Si_r}.$$
ii)Define the following symmetric bilinear form on
${\U}^-_n$ (the \textit{Green's scalar product}) :
$$\langle \f_\m,\f_{\m'}\rangle=v^{-2\; \mathrm{dim}\;\mathrm{Aut}(\m)} \frac{(1-v^2)^{|\m|}}{|\mathrm{Aut}(\m)|} \delta_{\m,\m'},$$
where $\mathrm{Aut}(\m)$ stands for
  the group of automorphism of any representation in the orbit
  $O_\m$ and $|\sum a_i^l [i;l)|=\sum_{i,l}la^l_i$.
It is natural to consider the restriction of this
scalar product $(\,,\,)$ on $\U^-_n$ to $\mathbf{R}\stackrel{i}{\simeq}
\Gamma$. Let $\mathcal{M}^{\mathrm{per}}$ denote the set of
multisegments of the form $\m=\sum a^l_i [i;l)$ such that $a^l_i=a^l_j$ for all $i,j$. By \cite{S}, Proposition 2.4 we have
$$\mathbf{R}=\big(\bigoplus_{\m \not\in \mathcal{M}^{\mathrm{per}}}
\A\b_\mathbf{m}\big)^\perp.$$
Hence the restriction of $(\,,\,)$ to $\mathbf{R}$ is nondegenerate. When $n=1$
this restriction coincides, up to a constant, with the Hall-Littlewood scalar
product. Let $(p_\mu)_{\mu \in \Pi}$ be the basis of power-sum symmetric functions and let $z_{(1^{m_1}2^{m_2}\cdots)}=\prod_i m_i!i^{m_i}.$
\begin{conj}The restriction of Green's scalar product on $\mathbf{R} \subset
\U^-_n$ is given by
$$(p_\lambda,p_\mu)=\delta_{\lambda,\mu}z_\lambda v^{-2(n-1)|\lambda|}
(1-v^2)^{n|\lambda|}\prod_{i=1}^{l(\lambda)}
\frac{1-v^{-2n\lambda_i}}{(1-v^{-2\lambda_i})^2}.$$
\end{conj}
This scalar product can be seen as a higher-rank analogue of the
Hall-Littlewood scalar product.
\section{Uglov's Fock spaces}
\paragraph{4.1} Let $n,l$ be positive integers and let $\mathbf{s}_l \in \Z^l$.
Following \cite{JMMO}, Uglov attached to this data an integrable $\U_v(\widehat{\mathfrak{sl}}_n)$-module $\Lambda^\infty_{\mathbf{s}_l}$ equipped with a distinguished
$\A$-basis $\{|\lambda_l,\mathbf{s}_l\rangle\}$, $\lambda_l \in \Pi^l$ (the
\textit{higher-level Fock space}, see \cite{U}, Section 1). The Fock space 
$\Lambda^\infty_{\mathbf{s}_l}$ is
also endowed with an action of a Heisenberg algebra $\mathcal{H}$ generated by
operators $B_m$, $m \in \Z^*$ (see \cite{U}, Sections 4.2, 4.3). Moreover, the
$\U_v(\widehat{\mathfrak{sl}}_n)$-action and the $\mathcal{H}$-action commute.
\paragraph{Remark.} When $l=1$, Uglov's Fock space
coincides with the Fock space $\Lambda^\infty$ introduced in \cite{KMS}.
\paragraph{4.2} We now extend the action of $\U_v^-(\widehat{\mathfrak{sl}}_n)$
on $\Lambda^\infty_{\mathbf{s}_l}$ to an action of $\U^-_n$. We follow the
method of Varagnolo-Vasserot \cite{VV}, Section 5. Let $\U^-_\infty$ be the 
Hall algebra of the quiver of type $A_\infty$. It is known that $\U^-_\infty
=\U^-_v(\mathfrak{sl}_\infty)$. Let $f_i$, $i \in \Z$ be the standard generator
corresponding to the vertex $i$.
\paragraph{}We associate to each $\lambda_l=(\lambda^{(1)},\ldots 
\lambda^{(l)}) \in \Pi^l$ an l-tuple of Young tableaux $(T_1,\ldots,T_l)$ such
that
\begin{enumerate}
\item[i)] $T_d$ is of shape $\lambda^{(d)}$ for $d=1,\ldots,l$,
\item[ii)] The $(i,j)$-box
of $T_d$ is filled with content $s_d+i-j$.
\end{enumerate}
If $\lambda_l$ and $\mu_l$ are two $l$-multipartitions such that $\gamma=
\mu_l\backslash \lambda_l$ corresponds to a box with content $k \in \Z$, we
say that $\gamma$ is an \textit{addable $k$-box} of $\lambda_l$ and a
\textit{removable $k$-box} of $\mu_l$. Let $\gamma,\gamma'$ be two addable
$k$-boxes of $\lambda_l$. We say that $\gamma < \gamma'$ if $\gamma$ and
$\gamma'$ belong to $T_d$ and $T_{d'}$ respectively and $d <d'$.\\
Let $\lambda_l,\mu_l \in \Pi^l$ be such that $\mu_l\backslash
\lambda_l$ is a $k$-box. Define
\begin{equation*}
\begin{split}
N^>(\mu_l,\lambda_l)=&\#\{addable\;k-boxes\;\gamma'\;of\;\lambda_l\;such\;
that\;\gamma'>\gamma\}\\
&-\#\{removable\;k-boxes\;\gamma'\;of\;\lambda_l\;such\;
that\;\gamma'>\gamma\}.
\end{split}
\end{equation*}
\begin{prop} The following endows $\Lambda_{\mathbf{s}_l}^\infty$ with a
structure of a $\U^-_\infty$-module :
$$f_k\cdot |\lambda_l,\mathbf{s}_l\rangle =\sum_{\mu_l}
v^{N^>(\mu_l,\lambda_l)}|\mu_l,\mathbf{s}_l\rangle$$
where the sum ranges over all $\mu_l$ for which $\mu_l\backslash \lambda_l$ is
a $k$-box.
\end{prop}
\noindent
\textit{Proof.} Straightforward. \qed
\paragraph{}Define operators $\mathbf{k}_k \in \mathrm{End}\;
(\Lambda^\infty_{\mathbf{s}_l})$, $k \in \Z$ by
$\mathbf{k}_k\cdot |\lambda_l,\mathbf{s}_l\rangle=v^{N_k(\lambda_l)}
|\lambda_l,\mathbf{s}_l\rangle$
where 
$$N_k(\lambda_l)=\#\{addable\;k-boxes\;of\;\lambda_l\}-
\#\{removable\;k-boxes\;of\;\lambda_l\}.$$
\paragraph{}Now let $d \in \N^{(\Z)}$ and set $\overline{d}=(\overline{d}_1,
\ldots,\overline{d}_n)$ where $\overline{d}_i=\sum_{j \equiv i\;
(\mathrm{mod}\;n)}
d_j$. Let $V$ be a $\Z$-graded $\mathbb{F}$-vector space of dimension $d$ and
let $\overline{V}$ be the $\Z/n\Z$-graded $\mathbb{F}$-vector space with
$\overline{V}_i=\bigoplus_{j \equiv i} V_j$. The collection of subspaces
$\overline{V}_{\leq i}=\bigoplus_{j \leq i} V_j$ defines a filtration of
$\overline{V}$ whose associated graded is $V$. Set
$$E_{\overline{V},V}=\{x \in E_{\overline{V}}\;|\;x(\overline{V}_{\geq i})
\subset \overline{V}_{\geq i+1}\;\mathrm{for\;all\;}i\}.$$
Let $p:E_{\overline{V},V} \to E_V$ be the projection onto the graded. Let
$j: E_{\overline{V},V} \subset E_V$ be the closed embedding. Following
\cite{VV}, define a map $\gamma_d :\U^-_n[\overline{d}] \to \U^-_\infty[d]$
by 
\begin{align*}
\gamma_{d|v=q^{-1}}:\;\C_{G_{\overline{V}}}(E_{\overline{V}}) &\to
\C_{G_V}(E_V)\\
f &\mapsto q^{-h(d)}p_!j^*(f)
\end{align*}
where $h(d)=\sum_{i<j, i \equiv j}d_i(d_{j+1}-d_j)$.

\paragraph{}For all $\lambda_l \in \Pi^l$ and $x \in \U^-_n$ we put
\begin{equation}\label{E:VV}
x \cdot |\lambda_l,\mathbf{s}_l\rangle=\sum_d \big(\gamma_d (x) \prod_{j<i,j
\equiv i} \mathbf{k}_{i}^{d_j}\big)\cdot |\lambda_l,\mathbf{s}_l\rangle.
\end{equation}
Then (see \cite{VV} Section 6.2, and \cite{Ariki})
\begin{prop} Formula (\ref{E:VV}) defines a representation 
$\Xi:\U^-_n \to \mathrm{End}\;(\Lambda^\infty_{\mathbf{s}_l})$ which extends
Uglov's action of $\U^-_v(\widehat{\mathfrak{sl}}_n)$.
\end{prop}
\paragraph{Remarks.} i) The number $h(d)$ has the following interpretation. Let
$\mathcal{F}_d$ be the variety of filtrations of $\overline{V}$ whose
associated graded is of dimension $d$. Then $\mathrm{dim}\;T^*\mathcal{F}_d=
\mathrm{dim}\;G_{\overline{V}} + h(d)$.\\
ii) The map $\gamma_d$ is ``upper triangular'' in the following sense. Let
$x \in E_V$ and define $r(x) \in E_{\overline{V}}$ by $r(x)_i=
\bigoplus_{j \equiv i} x_j$. Then $\gamma_d(\f_\m)(x) \neq 0 \Rightarrow r(x)
\in \overline{\mathcal{O}_\m}$.
\paragraph{4.3}Let $\mathcal{H}^- \subset \mathcal{H}$ denote the subalgebra
generated by $B_{-m}$, $m \in \N^*$. Define an algebra isomorphism $j: \Gamma
\stackrel{\sim}{\to}\mathcal{H}^-$ by setting $j(p_m)=B_{-m}$, where $p_m$ is
the power-sum symmetric function. Recall the canonical map $i:\mathbf{R}
\stackrel{\sim}{\to} \Gamma$ from Theorem 3.5.
\begin{lem}We have $\Xi_{|\mathbf{R}}=j \circ i$.\end{lem}
\noindent
\textit{Proof (sketch).} This is shown in a way similar to \cite{VV}. We first
consider the ``limit'' $\bigotimes^\infty$ of $\mathbf{T}_{n,r}$ when $r \to
\infty$ (see \cite{VV}, Section 10). Then $\Lambda^\infty_{\mathbf{s}_l}$ is
naturally embedded in a certain quotient of $\bigotimes^\infty$ (see \cite{U},
Section 3.3). In particular, the $\mathbf{U}^-_n$-action on $\mathbf{T}_{n,r}$
induces an action on $\bigotimes^\infty$ and on
$\Lambda^\infty_{\mathbf{s}_l}$. Let $\Xi'$ denote this last action. It
follows from Theorem 3.5 and \cite{U}, Section 4 that $\Xi'_{|\mathbf{R}}=j
\circ i$.
Finally, an easy extension to the higher-level Fock space of the computation in
\cite{VV}, Lemma 10.1 shows that $\Xi'=\Xi$.\qed

\section{Canonical bases of Fock spaces}
\paragraph{5.1} We keep the settings of the previous Section.
Uglov has defined a semilinear involution $a \mapsto \overline{a}$ on $\Lambda^\infty_{\mathbf{s}_l}$ (\cite{U}, Section 4.4) and two canonical bases $\{\b^{\pm}_{\lambda_l}\}_{\lambda_l \in \Pi^l}$
characterized by the following properties :
$$\overline{\b^{\pm}}_{\lambda_l} =\b^{\pm}_{\lambda_l},$$
$$\b^+_{\lambda_l} \in |\lambda_l\rangle +v\bigoplus_{\mu_l}\S |\mu_l\rangle ,\qquad \b^-_{\lambda_l} \in |\lambda_l\rangle +v^{-1}\bigoplus_{\mu_l} \overline{\S}|\mu_l\rangle.$$
He furthermore computed the transition matrices $[\b^{\pm}_{\lambda_l}:
|\mu_l,\mathbf{s}_l\rangle]$. In particular we have the following result.
\begin{theo}[\cite{U}, 3.26]
$$\b^-_{\lambda_l}=\sum_{\mu_l} \mathbf{P}^{-,\mathbf{s}_l}_{\mu_l,\lambda_l}
|\mu_l,\mathbf{s}_l\rangle.$$
\end{theo}
\paragraph{Remark.} When $l=1$, Uglov's canonical bases
coincide with the canonical bases considered by Leclerc-Thibon (\cite{LT}).
In that setting, the transition matrices above were first
obtained by Varagnolo and Vasserot \cite{VV}.
\paragraph{5.2}Let us now consider the nondegenerate scalar product $(\,,\,)$
on 
$\Lambda^\infty_{\mathbf{s}_l}$ for which $\{|\lambda_l,\mathbf{s}_l\rangle\}$
is orthonormal. Let $\{\b^{+*}_{\lambda_l}\}$ be the dual basis to $\{\b^+_{\lambda_l}\}$ with respect to the scalar product $(\,,\,)$.\\
\hbox to1em{\hfill}Define a semilinear isomorphism $\Lambda^\infty_{\mathbf{s}_l} \to \Lambda^\infty_{\mathbf{s}_l'},\;u \mapsto u'$ by $|\lambda_l,\mathbf{s}_l\rangle'=|\lambda_l',\mathbf{s}_l'\rangle$.
\begin{prop}[\cite{U}, 5.14] We have $(\b^{+*}_{\lambda_l})'=\b^-_{\lambda_l'}$.
\end{prop}
\paragraph{5.3} Let $\mathbf{B}_{\mathbf{s}_l}=\{\b^+_{\lambda_l}\}_{\lambda_l
\in \Pi^l}$ be the (positive) canonical basis of
$\Lambda^\infty_{\mathbf{s}_l}$.
\begin{theo} Let $\m \in \mathcal{M}$. Then $\b_\m \cdot|0,\mathbf{s}_l\rangle
\in \mathbf{B}_{\mathbf{s}_l} \cup \{0\}$.\end{theo}
\noindent
\textit{Proof.} Lemma 4.3 implies that the $\U^-_n$-action on
$\Lambda^\infty_{\mathbf{s}_l}$ is the same as that considered in \cite{S},
Section 4. The result follows from \cite{S}, Theorem 4.2. \qed
\paragraph{5.4}Define a map $\tau_{\mathbf{s}_l}: \mathcal{M} \to \Pi^l \cup
\{0\}$ by $\tau_{\mathbf{s}_l}(\m)=0$ if $\b_\m\cdot|0,\mathbf{s}_l\rangle =0$
and $\b_\m\cdot|0,\mathbf{s}_l\rangle=\b^+_{\tau_{\mathbf{s}_l}(\m)}$
otherwise.
This map is not easy to describe for a general $\mathbf{s}_l$. Nevertheless
we have the following result.
\paragraph{}Let $\m \in \mathcal{M}$ and let $\lambda=(\lambda^{(1)},\ldots,
\lambda^{(n)})$ be the associated $n$-multipartition. Let
$r=\sum_i l(\lambda^{(i)})$. Set $\mathbf{i}=(1^{l(\lambda^{(1)})},
2^{l(\lambda^{(2)})},\ldots) \in \mathcal{A}^n_r$ and
$$\mathbf{p}=(\lambda^{(1)}_1,\ldots,\lambda^{(1)}_{l(\lambda^{(1)})},
\lambda^{(2)}_1,\ldots)\in (\Z^+)^r.$$
Finally, let $\mathcal{M}_\mathbf{i}(\mathbf{p})=(p^{(1)},\ldots,p^{(l)})$ and
let $r_i \in \Z/n\Z$ be the residue of $p^{(i)}$.
\begin{lem} Suppose that $s_1 \gg s_2 \gg \cdots \gg s_l$ and that $s_i \equiv
r_i\;(\mathrm{mod}\;n)$ for $i=1,\ldots l$. Then $\tau_{\mathbf{s}_l}
(\m)=\mathcal{M}_\mathbf{i}(\mathbf{p})$.
\end{lem}
\noindent
\textit{Proof.} See appendix.\qed
\paragraph{5.5}\textit{Proof of Theorem 1.} Let $\sigma \in \Si_{\mathbf{i}}
\backslash \widehat{\Si}_r / \Si_{\mathbf{i}}$. It follows from Corollary 3.5
that $J^\mathbf{i}_{\lambda,\sigma} =J_{\lambda,\m}$ if there exists $\m \in
\mathcal{M}$ such that $\m^\mathbf{i}=\sigma$ and
$J^\mathbf{i}_{\lambda,\sigma}=0$ otherwise. From Section 3.4 we see that
$$\m^\mathbf{i} =\sigma \Longleftrightarrow \mathbf{i}\bullet \sigma \in (\Z^+)^r\;
\mathrm{and}\; \m=\sum_{j=1}^r [\mathbf{i}_j;(\mathbf{i}\bullet \sigma)_j).$$
Now we compute $J_{\lambda,\m}$. Let $l,\mathbf{p},(r_i)_{i=1}^l$ be
associated to $\m$ as in Section 5.4. Let $\mathbf{s}_l=(s_1,\ldots,s_l)$
be in the asymptotic region $s_1 \gg s_2 \gg \cdots \gg s_l$ and satisfy
$s_i \equiv r_i\;(\mathrm{mod}\;n)$ for all $i$. We evaluate both sides of
(\ref{E:0def}) on $|0,\mathbf{s}_l\rangle \in \Lambda^\infty_{\mathbf{s}_l}$.
On the one hand, it follows from Lemma 4.3 and Uglov's description of the
action of the Heisenberg algebra \cite{U}, Proposition 5.3 that
$$a_\lambda\cdot|0,\mathbf{s}_l\rangle=\sum_{\mu_1,\ldots,\mu_l} c^\lambda_{\mu_1,\ldots,\mu_l} v^{\sum (b-1)|\mu_b|} \bigg( \sum_{\nu_1,\ldots,\nu_l}
e_{\nu_1,\mu_1}\cdots e_{\nu_l,\mu_l} |(\nu_1,\ldots,\nu_l),\mathbf{s}_l\rangle
\bigg)$$
where $e_{\nu_i,\mu_i} \in \Z[v^{-1}]$ are defined by the relations
$s_{\mu_i}\vac=\sum_{\nu_i} e_{\nu_i,\mu_i}|\nu_i\rangle$ in the \textit{level
l=1} Fock space representation of $\U^-_n$. But by \cite{LT}, Theorem 6.9 we
have $s_{\mu_i}\cdot\vac=\b^-_{n\mu_i}$ and thus $e_{\nu_i,\mu_i}=\mathbf{P}^-_{\nu_i,n\mu_i}$. On the other hand, from Theorem 5.3 we have
$$\sum_\n J_{\lambda,\n} \b_\n \cdot|0,\mathbf{s}_l\rangle=
\sum_{\n, \tau_{\mathbf{s}_l}(\n) \neq 0} J_{\lambda,\n}\b^+_{\tau_{\mathbf{s}_l}(\n)}.$$
In particular, $J_{\lambda,\m}=\big((\b^+_{\tau_{\mathbf{s}_l}(\m)})^*,
a_\lambda |0,\mathbf{s}_l\rangle\big)$. But by Lemma 5.4 and Proposition~5.2,
$$(\b^+_{\tau_{\mathbf{s}_l}(\m)})^*=(\b^+_{\mathcal{M}_\mathbf{i}(p)})^*=
(\b^+_{\mathcal{M}(\sigma)})^*=(\b^-_{\mathcal{M}(\sigma)'})'.$$
Using the relations $(\overline{u},v)=(u',\overline{v'})$ for any $u,v \in
\Lambda^\infty_{\mathbf{s}_l}$ (\cite{U}, Proposition 5.13) and
$\overline{a_\lambda \cdot|0,\mathbf{s}_l\rangle}=a_\lambda \cdot|0,\mathbf{s}_l\rangle$
(\cite{U}, Proposition 4.2) we get
$$J_{\lambda,\m}=(\b^-_{\mathcal{M}(\sigma)'},a_\lambda \cdot|0,\mathbf{s}_l\rangle').$$
Now, from \cite{LT}, Theorem 7.13 i) we have $(\b^-_{n\mu})'=(-v)^{(n-1)|\mu|}
\b^-_{n\mu'}$ in the level $l=1$ Fock space. Thus
\begin{equation*}
\begin{split}
a_\lambda&\cdot|0,\mathbf{s}_l\rangle=\\
=&(-v)^{(n-1)|\lambda|}\sum_{\mu_1,\ldots,\mu_l}c^\lambda_{\mu_1,\ldots,\mu_l}
v^{\sum_{b=1}^l (b-1)|\mu_b|}\bigg(\sum_{\nu_1,\ldots,\nu_l} \mathbf{P}^-_{\nu_1, n\mu_1'}\cdots \mathbf{P}^-_{\nu_l,n\mu_l'}|\nu,\mathbf{s}_l\rangle \bigg)
\end{split}
\end{equation*}
where $\nu=(\nu_l, \ldots,\nu_1)$. The theorem follows. \qed
\section{On the center of $\U^-_n$}
\paragraph{}In this section we give a simple geometric characterization of the
central subalgebra $\mathbf{R} \subset \U^-_n$ in terms of the maps $\gamma_d:
\U^-_n \to \U^-_\infty$ defined in Section 4.2.
\paragraph{6.1} Let $d \in \N^{(\Z)}$ such that $d_i \in \{0,1\}$ for all $i$.
Then $d$ is the dimension of a unique (noncyclic) multisegment
 $\n_d=\sum_{k=1}^t [i_k;l_k)$ in $\Z$ satisfying the following
condition :
\begin{equation}\label{E:cond}
\forall\;j,k\qquad [i_k,l_k) \cup [i_j,l_j)\;is\;not\;a\;segment.
\end{equation}
Let $V_d$ be a $\Z$-graded $\F$-vector space of dimension $d$.
Set $l(d)=\sum_k (l_k-1)$.
Note that it follows from (\ref{E:cond}) that $E_{V_d}$ has a unique open
$G_{V_d}$-orbit, say $\mathcal{O}_d$.
\begin{lem}Suppose that $i_1 \gg i_2 \gg \cdots \gg i_t$ and set $\mathbf{i}_t
=(i_1,\ldots,i_t)$. Then for any $\f \in \U^-_\infty[d]$ we have
$$\f\cdot|0,\mathbf{i}_t\rangle=v^{-l(d)}\f_{|\mathcal{O}_d}|((l_1),\ldots,(l_k)),
\mathbf{i}_t\rangle.$$
\end{lem}
\noindent
\textit{Proof.} Note that $E_{V_d}=\prod_{k=1}^t E_{V_d(k)}$ where $V_d(k)=
\bigoplus_{l=0}^{l_k-1} \F V_{i_k+l}$. Let $f_k \in E_{V_d(k)}$ for $k=1,\ldots
,t$. From (\ref{E:cond}) and Section 4.2 we deduce that
$$f_1\cdots f_t \cdot|0,\mathbf{i}_t\rangle=\sum_{\nu_1,\ldots,\nu_t}
d_1(\nu_1)\cdots d_t(\nu_t) |(\nu_1,\ldots,\nu_t),\mathbf{i}_t\rangle$$
where $f_k|0,i_k\rangle=\sum_{\nu} d_k(\nu)|\nu,\mathbf{i}_k\rangle$ in the
level $l=1$ Fock space. But from \cite{VV}, Proposition 5., it is easy to
see that $f_k \cdot|0,\mathbf{i}_k
\rangle=v^{-(l_k-1)} (f_k)_{|\mathcal{O}_d(k)}|(l_k),i_k\rangle$ where
$\mathcal{O}_d(k) \subset E_{V_d(k)}$ is the open orbit.\qed
\paragraph{6.2} Recall the element $a_\lambda=i^{-1}(s_\lambda) \in
\mathbf{R}$. For any $\lambda,\mu \in \Pi$ let $K^\lambda_\mu \in \N$ be the
Kostka number. 
\begin{theo} Let $d \in \N^{(\Z)}$ such that $d_i \in \{0,1\}$. Then
$$\gamma_d(a_\lambda)_{|\mathcal{O}_d}=v^{l(d)+h(d)}
K^\lambda_{(u_1,\ldots,u_t)}$$
if there exists $i_k,u_k \in \Z$, $k=1,\ldots,t$ such that
$\n_d=\sum_{k=1}^t [i_k;nu_k)$, and $\gamma_d(a_\lambda)_{|\mathcal{O}_d}=0$
otherwise.
\end{theo}
\noindent
\textit{Proof.} Without loss of generality we may assume that
$\n=\sum_{k=1}^t [i_k;l_k)$ where $i_1>i_2>\cdots >i_t$.
Choose $d'=\cup_{k=1}^t [i'_k;l_k)$ where $i'_k \equiv i_k\;
(\mathrm{mod}\;n)$ and $i'_1 \gg i'_2 \gg \ldots \gg i'_t$. Let $\xi:
E_{V_d'} \stackrel{\sim}{\to} E_{V_d}$ be the obvious isomorphism. Then
$\xi \circ \gamma_{d'}=\gamma_d$. Now let us consider the Fock space
$\Lambda^\infty_{\mathbf{i}_t'}$ where $\mathbf{i}'_t=(i'_1,\ldots,i'_t)$.
Using \cite{U}, Proposition 5.3 we have
\begin{equation*}
\begin{split}
\big(a_\lambda\cdot|0,\mathbf{i}'_t\rangle,|((l_1),\ldots,(l_t)),&\mathbf{i}'_t
\rangle\big)\\
&=\sum_{\mu_1,\ldots,\mu_t} c^\lambda_{\mu_1,\ldots,\mu_t}
v^{\sum_b(b-1)|\mu_b|} \mathbf{P}^-_{(l_1),n\mu_1}\cdots \mathbf{P}^-_{(l_t),
n\mu_t}\\
&=\sum_{\mu_1,\ldots,\mu_t} \delta_{(l_1)=n\mu_1}\cdots \delta_{(l_t)=n\mu_t}
c^\lambda_{\mu_1,\ldots,\mu_t} v^{\sum_b(b-1)|\mu_b|}
\end{split}
\end{equation*}
Note that for any $u_1,\ldots,u_t \in \Z$ we have
$c^\lambda_{(u_1),\ldots,(u_t)}=K^\lambda_\mu$ where $\mu \in \Pi$
is the partition with parts $\{u_1,\ldots,u_t\}$.\\
\hbox to1em{\hfill}On the other hand, by Lemma 6.1
$$\big(a_\lambda\cdot|0,\mathbf{i}'_t\rangle,|((l_1),\ldots,(l_t)),
\mathbf{i}'_t \rangle\big)=v^{\epsilon(d',\mathbf{i}'_t)-l(d)}\gamma_{d'}
(a_\lambda)_{|\mathcal{O}_{d'}}$$
where
$$\epsilon(d',\mathbf{i}'_t)=
\sum_{l=1}^t\sum_{j \equiv i'_l;j <i'_l} d'_j.$$
The result now follows from the easily checked identity
$$\epsilon(d',\mathbf{i}'_t)=\sum_b(b-1)|\mu_b|+h(d')$$
when there exists $u_k \in \N$, $k=1, \ldots,t$ such that
$d'=\cup_k [i'_k,nu_k)$ and $\mu_k=(u_k)$.\qed
\paragraph{Remark.} It follows from Remark 4.2 ii) that the previous theorem
gives a characterization of the central element $a_\lambda$.
\section{Appendix}
\paragraph{}In this appendix we prove Lemma 5.4.
\paragraph{A.1} As in \cite{U}, Section 4, define a partial order on $\Pi^l$
(depending on $\mathbf{s}_l$) as follows. Let $\mu=(\mu^{(1)},\ldots,
\mu^{(l)}) \in
\Pi^l$. Set $k_i^{(d)}=\mu_i^{(d)}+s_d+1-i$ for $d=1,\ldots,l$ and $i \in \N$.
Let us write $k_i^{(d)}=c_i^{(d)} -nm_i^{(d)}$ where $c_i \in \{1,\ldots,n\}$,
and let $\mathbf{k}=(k_1 > k_2 > \cdots)$ be the ordered sequence whose
underlying set is $\{c_i^{(d)}+n(d-1)-nlm_i^{(d)}|\;i \in \N,\;d=1,\ldots,l\}$.
Let $s =s_1+\cdots+s_l$. It is easy to see that $k_i=s+1-i$ for $i \gg 0$ and
we denote by $\zeta(\mu)$ the partition such that $\zeta(\mu)_i=k_i-s+i-1$.
Now let $\mu,\nu \in \Pi^l$. By definition, we set $\mu \leq \nu$ if
$\zeta(\mu)\leq \zeta(\nu)$.
\paragraph{A.2}From now on we assume that $v=1$.
\paragraph{}It is more convenient to work with a different basis than
$\{\f_\n\}$. Let $\n \in \mathcal{M}$ and let $x \in \mathcal{O}_\n$. Set
$V_k =\mathrm{Ker}\;x^k$ and let $\alpha^1,\ldots,\alpha^r \in \N^{\Z/n\Z}$ be
such that
$$\mathrm{dim}\;V_k=\alpha^1+\cdots+\alpha^k,\qquad k=1,\ldots,r$$
and $\mathrm{dim}\;V_r=\mathrm{dim}\;\n$. Let $\f_{\alpha^i} \in \U^-_n$ be the
characteristic function of the trivial representation of the quiver
$\tilde{A}_{n-1}$ on $V_{\alpha^i} \simeq V_i/V_{i-1}$.
\begin{lemma}[\cite{VV}, Section 13] We have $\f_{\alpha^1}\cdots\f_{\alpha^r}
\in \f_{\n} + \bigoplus_{\l <\n}\N \f_{\l}$.
\end{lemma}
Now let $\n,\l \in \mathcal{M}$ such that $\mathrm{dim}\;\n=\mathrm{dim}\;\l$.
Let $(\beta^k)$ and $(\gamma^k)$ be the sequences of dimensions attached as
above to $\n$ and $\l$ respectively. If $\mathbf{u},\mathbf{v} \in \Z/n\Z$ we
write $\mathbf{u} \leq \mathbf{v}$ if $\mathbf{u}_i \leq \mathbf{v}_i$ for all
$i \in \Z/n\Z$.
\begin{lemma} We have $\n \geq \l$ if and only if 
\begin{equation}
\beta^1+\cdots+\beta^k\leq \gamma^1+\cdots+\gamma^k\;\qquad{for\;all\;}
k.\tag{a}
\end{equation}
\end{lemma}
\noindent
\textit{Proof.} Straightforward.\qed
\paragraph{}We will write $(\beta^k) \leq (\gamma^k)$ if (a) holds and if
$\sum_k \beta^k=\sum_k \gamma^k$. Let
$(\alpha^1,\ldots,\alpha^r)$ be the sequence attached to $\m$. We will
first prove
\begin{align}
\f_{\alpha^1}\cdots \f_{\alpha^r}\cdot|0,\mathbf{s}_l\rangle
&\in \N^*|\mathcal{M}_\mathbf{i}(\mathbf{p}),
\mathbf{s}_l\rangle
+ \bigoplus_{\mu \not\geq \mathcal{M}_\mathbf{i}(\mathbf{p})} \N|\mu,
\mathbf{s}_l\rangle
\tag{b}\\
\f_{\beta^1}\cdots \f_{\beta^r}\cdot|0,\mathbf{s}_l\rangle&\in
\bigoplus_{\mu \not>\mathcal{M}_\mathbf{i}
(\mathbf{p})} \N|\mu,\mathbf{s}_l\rangle\qquad {for\;all\;}(\beta^k) >
(\alpha^k).
\tag{c}
\end{align}

\begin{lemma} Let $\mu=(\mu^{(1)},\ldots,\mu^{(l)}) \in \Pi^l$ and let
$\beta \in N^{\Z/n\Z}$. We have
$$\f_\beta\cdot|\mu,\mathbf{s}_l\rangle=\sum_\nu |\nu,\mathbf{s}_l\rangle$$
where the sum ranges over all multipartitions $\nu=(\nu^{(1)},\ldots,
\nu^{(l)})$ such that
\begin{enumerate}
\item[i)] $\nu^{(i)}\backslash \mu^{(i)}$ is a skew diagram with at most one
box in each row,
\item[ii)] The number of boxes in $\cup_i \nu^{(i)}\backslash\mu^{(i)}$ with
content $j\;\mathrm{mod}\;n$ is $\beta_j$.
\end{enumerate}
\end{lemma}
\noindent
\textit{Proof.} Let ${d} \in \N^{(\Z)}$ such that ${d} \equiv
\beta\;(\mathrm{mod}\;n)$. Then $\gamma_{d|v=1}(\f_{\beta})=
\overset{\to}{\prod}_i f_i^{(d_i)}$, where $\overset{\to}{\prod}$
denotes the ordered product from $-\infty$ to $\infty$ (see \cite{VV},
Remark 6.1) and where $f_i^{(d_i)}$ is the divided power. Moreover, for any
$\sigma \in \Pi^l$,
$$f_i \cdot |\sigma,\mathbf{s}_l\rangle=\sum_\gamma |\gamma,\mathbf{s}_l\rangle
$$
where the sum ranges over all $\gamma \in \Pi^l$ such that $\gamma \backslash
\sigma$ is an $i$-box. The Lemma now follows from Section 4.2. \qed
\paragraph{} Finally, recall that $s_1 \gg s_2 \gg \cdots \gg s_l$. It is
clear from the definition that for $\mu,\lambda \in \Pi^l$,
\begin{equation}\label{E:Aun}
\mu \geq \lambda \Rightarrow \;\exists\;k\;{such\;that}\;\mu^{(i)}=
\lambda^{(i)}\;{for}\;i=1,\ldots, k-1\;{and}\;\mu^{(k)} \geq
\lambda^{(k)}.\tag{d}
\end{equation}
Note that $\alpha^k_i$ is equal to the number of boxes with content $i$ in
the slice $s_k$ of the diagram $D_\mathbf{p}$ associated to $\mathbf{p}$.
Statements (b) and (c) now easily follow by Lemma 3 and by construction of
$\mathcal{M}_\mathbf{i}(\mathbf{p})$.
\paragraph{A.3} By \cite{U}, Theorem 2.4 it possible to
choose $s_l \gg s_{l+1} \gg \cdots \gg s_t$ for some $t \gg
0$ in such a way that $\b_\m |0,\mathbf{s}_t\rangle \neq 0$, where
$\mathbf{s}_t=(s_1,\ldots,s_t)$.
\begin{lemma} We have $\b_\m|0,\mathbf{s}_t\rangle=
\b^+_{\widetilde{\mathcal{M}}_\mathbf{i}(\mathbf{p})}$, where
$\widetilde{\mathcal{M}}_\mathbf{i}(\mathbf{p})
=\big(\mathcal{M}_\mathbf{i}(\mathbf{p}),0^{t-l}\big)$.
\end{lemma}
\noindent
\textit{Proof.} By Lemma 1, we have
$$\f_{\alpha^1}\cdots\f_{\alpha^r}\cdot|0,\mathbf{s}_t\rangle\in |\tau_{\mathbf{s}_t}(\m),\mathbf{s}_t\rangle
+
\bigoplus_{\mu < \tau_{\mathbf{s}_t}(\m)} \Z |\mu,\mathbf{s}_t\rangle.$$
But from (b) and (\ref{E:Aun}) it is clear that
$$\f_{\alpha^1}\cdots\f_{\alpha^r}\cdot|0,\mathbf{s}_t\rangle\in 
|\widetilde{\mathcal{M}}_\mathbf{i}(\mathbf{p}),\mathbf{s}_t\rangle
+
\bigoplus_{\mu \not\geq \widetilde{\mathcal{M}}_\mathbf{i}(\mathbf{p})}
 \Z |\mu,\mathbf{s}_t\rangle.$$
Hence $\tau_{\mathbf{s}_t}(\m)=
\widetilde{\mathcal{M}}_{\mathbf{i}}(\mathbf{p})$. \qed
\paragraph{}In particular,
$$\b_{\m} \cdot|0,\mathbf{s}_t\rangle \in 
|\widetilde{\mathcal{M}}_\mathbf{i}(\mathbf{p}),\mathbf{s}_t\rangle
+ \bigoplus_{\mu <\widetilde{\mathcal{M}}_\mathbf{i}(\mathbf{p})} \N |\mu,
\mathbf{s}_l\rangle.$$
Consider the projection $\pi: \Lambda^\infty_{\mathbf{s}_t} \to
\Lambda^\infty_{\mathbf{s}_l}$ given by
$$|\big(\mu^{(1)},\ldots,\mu^{(t)}\big),\mathbf{s}_t\rangle \mapsto
\cases |\big(\mu^{(1)},\ldots,\mu^{(l)}\big),\mathbf{s}_l\rangle &
if\;\mu^{(j)}=0\;for\;j>l\\
0 & otherwise \endcases$$
It is clear from (\ref{E:VV}) that $\pi(\b_\m\cdot|0,\mathbf{s}_t\rangle)=
\b_\m \cdot|0,\mathbf{s}_l\rangle$. Hence
$$\b_\m\cdot|0,\mathbf{s}_l\rangle \in |\mathcal{M}_\mathbf{i}(\mathbf{p}),
\mathbf{s}_l\rangle + \bigoplus_{\mu < \mathcal{M}_\mathbf{i}(\mathbf{p})}
\N |\mu,\mathbf{s}_l\rangle.$$
This proves Lemma 5.4 \qed

\vspace{1cm}
\centerline{\textbf{Acknowledgements}}
I am grateful to J. Dat for interesting discussions and for pointing out
to me the work \cite{Hai1}, \cite{Hai2}, to B. Leclerc, A. Schilling and
E. Vasserot for valuable advice at various stages of this work.
I would like to thank the MIT mathematics department for its hospitality.
This research was partially conducted for the Clay Mathematics Institute.

\small{}
\vspace{4mm}
Olivier Schiffmann, MIT, 77 Massachusetts Avenue, CAMBRIDGE 02139, USA;\\
\hbox to11em{\hfill}email:\;\texttt{schiffma@math.mit.edu}

\end{document}